\documentclass{amsart}
\usepackage{amssymb}
\input xypic
\xyoption{all}
\renewcommand{\O}{\mathcal{O}}
\newcommand{\ord}{\operatorname{ord}}
\newcommand{\isom}{\cong}
\newcommand{\lsplt}{\Theta}
\newcommand{\XX}{\mathcal{X}}
\newcommand{\hr}{H_{\textup{rig}}}
\newcommand{\dr}{\textup{dR}}
\newcommand{\et}{{\textup{\'et}}}
\newcommand{\Q}{\mathbb{Q}}
\newcommand{\C}{\mathbb{C}}
\newcommand{\cy}{\mathcal{Y}}
\newcommand{\cz}{\mathcal{Z}}
\newcommand{\tA}{\tilde{A}}
\newcommand{\tB}{\tilde{B}}

\def\het^#1(#2,#3){H_\et^#1(#2,\Q_p(#3))}
\def\hset_#1^#2(#3,#4){H_{#1,\et}^#2(#3,\Q_p(#4))}
\newcommand{\res}{\operatorname{Res}}
\newcommand{\inject}{\hookrightarrow}
\newcommand{\pair}[1]{\langle #1 \rangle}
\newcommand{\hdr}{H_{\dr}}
\newcommand{\hcr}{H_{\textup{cr}}}
\newcommand{\gm}{{\mathbb{G}_{\textup{m}}}}
\newcommand{\ddr}{D_{\dr}}
\newcommand{\dcr}{D_{\textup{cr}}}
\newcommand{\hf}{H_{\textup{f}}}
\newcommand{\nekovar}{Nekov\'a\v r}
\newcommand{\Nek}{\textup{Nek}}
\newcommand{\CG}{\textup{CG}}
\newcommand{\poincare}{Poincar\'e}
\newcommand{\uni}{\varpi}
\newcommand{\lunip}{\log_{\uni}^\prime}

\newcommand{\bsY}{\overline{\cy}}
\newcommand{\bsZ}{\overline{\cz}}
\newcommand{\etale}{\'etale}
\newcommand{\Xbar}{\overline{X}}
\newcommand{\Hom}{\operatorname{Hom}}
\newcommand{\spec}{\operatorname{Spec}}
\newcommand{\Gal}{\operatorname{Gal}}
\newcommand{\ndiv}{\nmid}
\newcommand{\Div}{\operatorname{Div}}
\newcommand{\Olog}{\Omega_{\log}^1} 
\newcommand{\lag}{\Psi}

\newtheorem{theorem}{Theorem}[section]

\newtheorem{proposition}[theorem]{Proposition}
\newtheorem{lemma}[theorem]{Lemma}
\newtheorem{corollary}[theorem]{Corollary}
\theoremstyle{definition}
\newtheorem{definition}[theorem]{Definition}

\numberwithin{equation}{section}
\begin{document}
\title{The $p$-adic height pairings of Coleman-Gross and of Nekov\'a\v r}
\author{Amnon Besser}
\maketitle

\section{introduction}
\label{sec:intro}

In~\cite{Col-Gro89}, Coleman and Gross proposed a definition of
a $p$-adic height pairing on curves over number fields with good
reduction at primes above $p$. The pairing was defined as a sum of
local terms and the most interesting terms are the ones corresponding
to primes above $p$ where the definition depends on Coleman's theory
of $p$-adic integration. Later, \nekovar\ constructed in~\cite{Nek93} a
general $p$-adic height pairing for Galois representations which are
cristalline at primes above $p$, and satisfying certain additional
technical assumptions.

Wintenberger raised the question of the equality of the two height
pairings for curves. In other words: does the \nekovar\ height pairing
applied to $H_\et^1$ of a curve recover the Coleman-Gross pairing. The
purpose of this small note is to answer this in the
affirmative. More precisely, let $F$ be a number field and let $X$ be
a smooth and proper curve over $F$, with good reduction at all places $v$
above a fixed prime $p$. To define either the height of Coleman and
Gross or the height of \nekovar\ the following choices must be made.

\begin{itemize}
\item  A ``global log''- a continuous idele class character
 \begin{equation*}
   \ell: \mathbb{A}_F^\times/F^\times \to \Q_p\;.
 \end{equation*}
\item for each $v|p$ a choice of a subspace $W_v\in \hdr^1(X\otimes
  F_v/F_v)$ complementary to the space of holomorphic forms.
\end{itemize}
For the definition of the Coleman-Gross height we must insist that the
local characters $\ell_v$, for $v|p$, are ramified in the sense that
they do not vanish on the units in $F_v$ (this seems to be overlooked
in~\cite{Col-Gro89}).

With these choices the Coleman-Gross height pairing is a pairing
\begin{equation*}
  h_{\CG}: \Div_0(X) \times \Div_0(X) \to \Q_p\;,
\end{equation*}
where $\Div_0(X)$ denotes the group of zero divisors on $X$ (defined
over $F$),
while the \nekovar\ height pairing is a pairing
\begin{equation*}
  h_{\Nek}: \hf^1(F,V) \times \hf^1(F,V) \to \Q_p\;,
\end{equation*}
where $V=V_p(J_X)\isom \het^1(X,1)$ is the Tate-module of the
jacobian of $X$ and $\hf$ is finite cohomology in the sense of Bloch
and Kato. The relation between the two pairings is provided by
the \etale\ Abel-Jacobi map
\begin{equation*}
  \alpha_X: \Div_0(X) \to \hf^1(F,V)\;.
\end{equation*}
Our main result is the following.
\begin{theorem}\label{mainthm}
  With the same choices of $(\ell,W_v)$ we have
  \begin{equation*}
    h_{\Nek}(\alpha_X(y),\alpha_X(z))=h_{\CG}(y,z)\;.
  \end{equation*}
\end{theorem}

The proof goes as follows. We first recall the construction of the
Coleman-Gross height in Section~\ref{sec:colgros}. This is expressed
in terms of local heights. Likewise, the \nekovar\ height can be
expressed as a sum of local heights, although it also has a global
description which we will not need here. The decomposition into local
heights depends on a choice of a ``mixed extension''. We recall this
in Section~\ref{sec:mixed} as well as the particular mixed extension we will
use. With this mixed extension, it is known that the local heights at
places not dividing $p$ are equal to the local heights of
Coleman-Gross. This leaves the comparison of the local heights above
$p$. This is done in the rest of the paper.

We would like to thank Wintenberger and \nekovar\ for suggesting this
problem to us on various occasions.

\section{The Coleman-Gross height pairing}
\label{sec:colgros}

We first recall the theory of p-adic height pairings due to
Coleman and Gross \cite{Col-Gro89}.

Recall that we have the character $\ell: \mathbb{A}_F^\times/F^\times \to \Q_p$.
 One deduces from $\ell$ the following data:
 \begin{itemize}
 \item For any place $v\ndiv p$ we have $\ell_v(\O_{F_v}^\times)=0$ for
   continuity reasons, which implies that $\ell_v$ is completely
   determined by the number $\ell_v(\pi_v)$, where $\pi_v$ is any
   uniformizer in $F_v$.
 \item For any place $v|p$ we can decompose
   \begin{equation*}
     \xymatrix{
     {\O_{F_v}^\times}  \ar[rr]^{\ell_v} \ar[dr]^{\log_v} & &   \Q_p\\
         & F_v\ar[ur]^{t_v}
           }
   \end{equation*}
   where $t_v$ is a $\Q_p$-linear map. Since we assume that $\ell_v$
   is ramified it is then possible to extend
   $\log_v$ to $F_v^\times$ in such a way that the diagram remain
   commutative.
 \end{itemize}

The height pairing is a sum of local terms, $h_{\CG}(y,z)= \sum_v
h_{\CG}^v(y,z)$ over all finite places
$v$. When $v\ndiv p$ the local term is given~\cite[(1.3)]{Col-Gro89} by
\begin{equation}
  \label{locCG}
  h_{\CG}^v(y,z) = \ell_v(\pi_v) \cdot (y,z)
\end{equation}
where $(y,z)$ denotes the intersection multiplicity of the extension
of $y$ and $z$ to a regular model of $X$ over $\O_{F_v}$.

We now describe the local contribution at a place $v|p$. We modify the
definitions of Coleman and Gross to work with integration theory over
$\overline{\Q}_p$ instead of over $\C_p$. Let $X$ be a
curve over $K=F_v$ with good 
reduction. Let $\Olog(X)$ be the space of one forms on $X$ (defined
over $K$) with at
most logarithmic singularities. For $\omega\in \Olog(X)$ we let
$\res(\omega)$ be the residue divisor. This is a divisor of degree $0$
defined over $K$ such that over
$\overline{K}$ it becomes $\sum_{P\in X} \res_P(\omega) P$. Let
$\Omega^1(X)$ be the space of global holomorphic 
forms on $X$. We note that if $\omega_1$ and $\omega_2\in \Olog(X)$
and $\res(\omega_1)=\res(\omega_2)$ then $\omega_1-\omega_2\in
\Omega^1(X)$.
There is a canonical projection 
\begin{equation}
  \lag:\Olog(X)\to \hdr^1(X)\;.\label{eq:lag}
\end{equation}
 In
\cite{Col-Gro89} it is given by certain requirements. Following
\cite{Bes00} we can describe it as follows:
\begin{proposition}\label{proplag}
  Let $\omega\in \Olog(X)$ and let $U\subset X$
  be a wide open subspace on which $\omega$ is holomorphic.
  The class $\lag(\omega)$ is the
  image of $\omega|_U \in \hdr^1(U)$ under the unique Frobenius
  equivariant section to
  the map $\hdr^1(X) \to \hdr^1(U)$.  
\end{proposition}

Now recall that we have at our disposal the complementary subspace $W=W_v$.
\begin{definition}\label{omd}
  For any divisor $y$ of degree $0$ on $X$ we let $\omega_W(y)\in
  \Olog(X)$ be the unique form satisfying $\res(\omega_W(y))=y$ and
  $\lag(\omega_W(y))\in W$.
\end{definition}

Now we need to use the theory of $p$-adic
integration. In~\cite{Col-Gro89} the theory of Coleman~\cite{Col-de88}
is used.
Several other $p$-adic integration theories have been developed in recent
years~\cite{Zar96,Colm96,Bes99,Vol01,Bes00}. In most of these theories
one can at least integrate algebraic differentials on all smooth algebraic
varieties over a
finite extension of $\Q_p$ and in this domain of integration they
are conjectured to be identical. In~\cite{Bes00} it is shown that the
theory developed there using~\cite{Vol01}, as well as the theory
of~\cite{Colm96}, coincide for curves with Coleman's integration
theory. For our needs these theories are equivalent.

The integration theory depends on a branch of the $p$-adic log. For
this we take the branch $\log_v$ defined by $\ell_v$ and we can extend
it in a unique way to $\overline{K}$. The properties of $p$-adic
integration, which are valid for all versions
except~\cite{Bes99,Col-de88} (which work in the rigid setting but
require good or close to good reduction), can be summarized as follows.
\begin{proposition}\label{funct}
  For any smooth algebraic variety $X/K$ and a holomorphic differential
  $\omega$ on $X$ defined over $K$ the theory associates a function
  $\int \omega: X(\overline{K})\to \overline{K}$ unique up to an
  additive constant, such that the following properties are satisfied.
  \begin{itemize}
  \item The function $\int \omega$ is locally analytic and satisfies
    $d \int \omega = \omega$.
  \item The integration is functorial with respect to arbitrary
    morphisms in the sense that if $f:Y\to X$ is defined over $K$, then
    we have $\int f^\ast \omega = f^\ast \int \omega$ up to a constant.
  \item if $\sigma\in \Gal(\overline{K}/K)$ and $x,y \in
  X(\overline{K})$, then $\int_{\sigma(x)}^{\sigma(y)} \omega =
  \sigma(\int_x^y \omega)$.
  \end{itemize}
\end{proposition}
The following easy corollary of the above properties is well known
(and is in fact used by Zarhin and Colmez to define the integral)
\begin{corollary}\label{funct1}
  Let $J$ be a commutative algebraic group over $K$, let $\omega$ be
  an invariant differential on $J$ and let $\omega_0$ be its value at
  the identity element. Let $\log: J(K) \to T_0(J)$ be the logarithm
  for the $p$-adic Lie group $J$, depending on the choice of a branch
  of the logarithm, taking values in the tangent space at $0$ to
  $J$. Then we have $\int_0^P \omega=\pair{\omega_0,\log(P)}$, where
  $\pair{~,~}$ is the duality between the cotangent and tangent space
  at the identity.
\end{corollary}

Suppose again that $X$ is a curve as before and let $\omega\in
\Olog(X)$. Then we obtain an integral $\int
\omega$ defined on $X(\overline{K})$ minus the singular points of the
form. The
sum of the values of $\int \omega$ on $z\in \Div_0(X)$, disjoint from
the singular points, is
well defined independent of the constant of integration and denoted
$\int_z\omega$. Furthermore, since $z$ and $\omega$ are
defined over $K$ we have $\int_z\omega \in K$.

\begin{definition}\label{locCol}
  The local height pairing at $v|p$ of Coleman and Gross is defined as
  follows: 
  Let $y$ and $z$ be two divisors of degree $0$ on $X$ with disjoint
  supports. Then their
  pairing is given by $h_{\CG}^v(y,z):=t_v(\int_{z}\omega_{W_v}(y))$.
\end{definition}

\section{The \nekovar\ height in terms of mixed extensions}
\label{sec:mixed}

\nekovar\ gives two expressions for the height pairing. The one which
is relevant for us is given in terms of mixed extensions and is
expressed as a sum of local terms. As we saw in the previous section
the Coleman-Gross height pairing is
also given as a sum of local terms. That the terms at places not
dividing $p$ are the same is well-known and we recall this here.

The height pairing of \nekovar\ is a map
\begin{equation*}
  h_{\Nek}: \hf^1(F,V) \times \hf^1(F,V^\ast(1)) \to \Q_p\;,
\end{equation*}
where $V$ is a continuous $G=\Gal(\overline{F}/F)$-representation
satisfying certain
conditions (see \cite[2.1.2]{Nek93}) and $\hf$ is finite cohomology in
the sense of
Bloch and Kato~\cite{Blo-Kat90}. To describe the height pairing we
assume that we are given extensions
of continuous $G$-representations,
\begin{equation}\label{shortA}
  0\to V \to E_1 \to A \to 0
\end{equation}
and
\begin{equation}\label{shortB}
  0\to B(1) \to E_2 \to V \to 0\;.
\end{equation}
Here, $A$ and $B$ are finite dimensional $\Q_p$-representations of
$G$ which are trivialized by a finite extension (this is more general
then in \cite{Nek93} and is inspired by~\cite{Scho94}). Suppose for
now that $G$ acts trivially on $A$ and $B$. Then, such sequences
yield (the second by first 
dualizing) elements $[E_1]\in \Hom(A,H^1(F,V))$ and $[E_2]\in B\otimes
H^1(F,V^\ast(1))$. We assume that these in fact belong to
$\Hom(A,\hf^1(F,V)$ and $B\otimes \hf^1(F,V^\ast(1))$
respectively. The height pairing produces out of these two extensions
an element $h(E_1,E_2)\in \Hom(A,B)$. The construction is functorial
with respect
to maps $A'\to A$ and $B\to B'$. To get the height pairing itself one
restricts to the case $A=B=\Q_p$. The more general setup is convenient
for example in geometric situations, as we will see.

To describe the height pairing one chooses a mixed extension. This is
an embedding of the sequences \eqref{shortA} and \eqref{shortB} inside
a commutative diagram with exact rows and columns as follows
\cite[p. 159]{Nek93}:
\begin{equation}\label{mixext}
\xymatrix{
 & & 0 \ar[d] & 0 \ar[d] & \\
0 \ar[r]& B(1)\ar[r] \ar@2{-}[d] & E_2 \ar[r]^{\pi} \ar[d] & V\ar[r] \ar[d] & 0\\
0 \ar[r]& B(1)\ar[r]  & E \ar[r] \ar[d] & E_1 \ar[r] \ar[d] & 0\\
 &  & A \ar@2{-}[r] \ar[d] & A \ar[d] & \\
 &  & 0  & 0  & 
}
\end{equation}

Having fixed such a mixed extension $E$ we obtain local mixed extensions
$E_v$ by restricting to the decomposition groups at the various places
$v$. Each of these local mixed extensions gives a local contribution
$h_{E_v}$ 
and the sum $h_E = \sum_v h_{E_v}$ is shown to depend only on $E_1$
and $E_2$ so we can define $h(E_1,E_2)=h_E$.
We will only need to describe the local contribution at
places above $p$, which we will do in the next section.

Suppose now that $X$ is a curve as in the statement of the theorem and
that $y$ and $z$ are two zero divisors on $X$ with disjoint
supports. Let $\cy$ and $\cz$ be the supports of $y$ and $z$
respectively. We will use an overline to denote extension of scalars
to the algebraic closure of the corresponding field. Then we obtain our mixed
extension as follows (compare~\cite[5.6]{Nek93}. Let
$\tA=(\Q_p^{\overline{\cy}})_0$ be the subspace of the space
of functions from $\overline{\cy}$ to $\Q_p$ where the sum of values
is $0$, and
let $\tB=(\Q_p^{\overline{\cz}})^0$ be the quotient of the
corresponding space for $\cz$ by the
subspace of constant functions.
\begin{definition}\label{basicmixed}
  The semi geometric mixed extension associated with the cycles $y$ and $z$
  is the mixed extension in \eqref{mixext} with 
  \begin{align*}
    V&=\het^1(\Xbar,1),\quad E_1=\het^1(\Xbar - \overline{\cy},1),\quad
    E_2=\het^1(\Xbar; \overline{\cz},1)\\
    \intertext{(\etale\ cohomology of
      $\Xbar$ relative to $\overline{\cz}$), and finally}\\
    E&=\het^1(\Xbar- \overline{\cy}; \overline{\cz},1)\;.
\end{align*}
 The associated geometric mixed
  extension is the extension obtained from the semi geometric mixed
  extension by the maps $A=\tA^G \inject \tA$ and $\tB \twoheadrightarrow
  B= \tB_G$.
\end{definition}
In this
mixed extension, $\tA$ is identified with the kernel of
$\hset_{\bsY}^2(\Xbar,1) \to \het^2(\Xbar,1)$ and $\tB$ with the
cokernel of $\het^0(\Xbar,0)\to \het^0(\Xbar;\bsZ,0)$.

Consider now $y,z\in \Div_0(X)$ as $y\in A$ and $z\in B^\ast$
($\Q_p$-dual) in the obvious manner. Then it is well known that
$[E_1](y)= \alpha_X(y)$ and $((z\otimes 1)\circ
[E_2])(1)=\alpha_X(z)$. We thus find 
\begin{equation*}
  h_{\Nek}(y,z)= z(h_E(y))=\sum_v z(h_{E_v}(y))\;.
\end{equation*}
\begin{proposition}
  For $v\ndiv p$ we have $z(h_{E_v}(y))=h_{\CG}^v(y,z)$, where $E$ is
  the mixed extension 
  described above and $h_{\CG}^v$ is the local height pairing of Coleman
  and Gross described in \eqref{locCG}.
\end{proposition}
\begin{proof}
  This is well known. See for example~\cite[Proposition 2.16]{Nek95}
  where a more general result is proved.
\end{proof}

To prove Theorem~\ref{mainthm}
it thus remains to prove the following result.
\begin{theorem}\label{locthm}
  For $v| p$ we have, with the same choice of the space $W_v$,
  $z(h_{E_v}(y))=h_{\CG}^v(y,z)$.
\end{theorem}

The proof of this theorem will occupy the rest of this paper. Before
proving this we want to make the following comment: Since we will
dualize several times with respect to cup products, one has to check
the signs very carefully here. We did not do this, hence, apriori,
Theorem~\ref{mainthm} is only true up to a global sign. This sign
can not be $-1$, however, as can be seen from the fact that both the
\nekovar\ height pairing and the Coleman-Gross height pairing have the
property that they vanish if either $y$ or $z$ are principal divisors.

\section{The local height at a place above $p$}

We now begin to compute the local
contribution to the height pairing at a place $v|p$ coming from the
mixed extension described in the previous section. In this section we
do this using a certain splitting (the map $w$ below), which will be
described in Section~\ref{sec:end}.

We begin by
recalling the general construction of \nekovar\ in
\cite[4.7]{Nek93}. We make several modifications to bring it to a form
suitable for comparison with the Coleman-Gross construction.

Let $K=F_v$. The Kummer map,
\begin{equation}\label{eq:kummer}
 K^\times \otimes \Q_p \xrightarrow{\alpha_K} H^1(K,\Q_p(1))\;,
\end{equation}
is an isomorphism and identifies
$\hf^1(K,\Q_p(1))$ with $\O_K^\times \otimes \Q_p$, so that
there exists a short exact sequence 
\begin{equation}
  \label{eq:tateshort}
  0\to \hf^1(K,\Q_p(1))\to H^1(K,\Q_p(1)) \xrightarrow{\ord} \Q_p \to 0\;.
\end{equation}
The $p$-adic logarithm $\log: \O_K^\times \to K$
defines an isomorphism
$\hf^1(K,\Q_p(1)) \isom K$. The branches of the $p$-adic
logarithms, i.e., extensions of $\log$ to $K^\times$, are therefore in
one to one correspondence with splittings of \eqref{eq:tateshort}. We
will denote the splitting induced by a branch $\log_\uni$ by
$\lunip$. We can view the local component of the idele class
character $\ell_v$ as a map $\ell_v: H^1(K,\Q_p(1)\to \Q_p$. Then,
the character $\ell_v$ factors as $t_v\circ \log_\uni$ if and only if
we have
$\ell_v = \ell_v\circ \lunip$. Note that in this case we have
\begin{equation}
  \label{whertvcomes}
  (K\to \hf^1(K,\Q_p(1)) \xrightarrow{\ell_v} \Q_p) =t_v\;.
\end{equation}
From now onward we assume that such a choice of
$\log_\uni$ has been made.

We now consider the diagram \eqref{mixext} as a
diagram of representations of $\Gal(\overline{K}/K)$ and assume that
$A$ and $B$ are trivial representations. Recall that  $\pi:
E_2 \to V$. We then have the following diagram
\begin{equation}\label{bigdiag}
\xymatrix{
&0\ar[d]&0\ar[d]&& \\
0\ar[r]&B\otimes \hf^1(K,\Q_p(1))\ar[r]^{j}\ar[d]^u&\hf^1(K,E_2)\ar@<1ex>@{.>}[l]^w
\ar[r]^\pi \ar[d]^i&\hf^1(K,V)\ar@2{-}[d]\ar[r]&0\\
0\ar[r]&B\otimes H^1(K,\Q_p(1))\ar[r]^j \ar[d]^{1\otimes
  \ord}\ar@<1ex>@{.>}[u]^{1\otimes\lunip}\ar[ddl]^{1\otimes \ell_v}
&\pi^{-1}(\hf^1(K,V))\ar@<1ex>@{.>}[l]^{w'}\ar@<1ex>@{.>}[u]^{\lsplt}
\ar[r]^\pi \ar[d]^\phi&\hf^1(K,V)\ar[r]&0\\
&B\ar@2{-}[r]\ar[d]&B\ar[d]&&\\
B&0&0&&
}
\end{equation}

This diagram is mostly copied from~\cite[4.8]{Nek93}.
The dotted lines in the diagram denote various splittings. We have
already discussed $\lunip$. The map $w$ depends on the choice of the
complementary subspace $W_v$ from the introduction (it
does not depend on a choice of $\log_\uni$) and
will be discussed later. The two other splittings are induced, $w'$ by
$w$ and $\lsplt$ by $\lunip$ in the following formal way.
\begin{definition}\label{splitdef}
  Let $x\in \pi^{-1}(\hf^1(K,V))$ and choose $y\in \hf^1(K,E_2)$ such
  that $\pi(y)=\pi(x)$. Then $\pi(x-i(y))=0$ so $x-i(y)=j(z)$ with
  $z\in B\otimes H^1(K,\Q_p(1))$. Then we define
  \begin{align*}
    w'(x)&=u(w(y))+z\\
    \intertext{and}
    \lsplt(x)&= y+ j\circ (1\otimes \lunip)(z)
  \end{align*}
\end{definition}
The following is an easy diagram chase away.
\begin{lemma}\label{4.2}
  The maps $w'$ and $\lsplt$ are well defined and constitute
  splittings of the respective diagrams. Furthermore, in diagram
  \eqref{bigdiag} the two possible ways of connecting $\hf^1(K,E_2)$
  with $B\otimes H^1(K,\Q_p(1))$ (in both directions) give the same
  map.
\end{lemma}

In~\cite[4.8]{Nek93} \nekovar\ describes the local height pairing
separately in the two cases $\ell_v$ ramified or not. However it is
easy to see that his two definitions coincide with the following.
\begin{definition}
  The local height pairing corresponding to the mixed extension $E$ is
  the element of $h_E\in \Hom(A,B)$ obtained as follows: We have $[E]\in
  \Hom(A,H^1(K,E_2))$, but in fact lies in
  $\Hom(A,\pi^{-1}(\hf^1(K,V))$. Then $h_E$ is the composition
  \begin{equation*}
    A\xrightarrow{[E]} \pi^{-1}(\hf^1(K,V)) \xrightarrow{w'} B\otimes
    H^1(K,\Q_p(1)) \xrightarrow{1\otimes \ell_v} B\;.
  \end{equation*}
\end{definition}

Since we are assuming the existence of a $\log_\uni$ through which
$\ell_v$ factors, we can write this in a different way.
\begin{lemma}
  We have $(1\otimes \ell_v) \circ w' = (1\otimes \ell_v) \circ u\circ
  w \circ \lsplt$.
\end{lemma}
\begin{proof}
With $x$, $y$ and $z$ as in Definition~\ref{splitdef} we have
\begin{equation*}
 u\circ  w \circ \lsplt(x)=   u(w(y+ (1\otimes \lunip)(z))) = u(w(y)+
 (1\otimes \lunip)(z))
\end{equation*}
\newcommand{\id}{\operatorname{Id}}
so that $u\circ  w \circ \lsplt(x)-w'(x)=  u((1\otimes
\lunip)(z))-z$. But the compatibility of $\log_\uni$ with $\ell_v$
precisely means that $\ell_v$ vanishes on the image of $u\circ \lunip
- \id$.
\end{proof}

As a consequence we have obtained the description of the local height pairing
as the composition
\begin{equation}\label{bigfactor}
  h_E = (1\otimes \ell_v) \circ u\circ w \circ
  \lsplt \circ [E]\;.
\end{equation}

\section{The connection with $p$-adic integration}
\label{sec:abel-jacobi}

In this section we go back to the geometric mixed extension of
Definition \ref{basicmixed} and
we show that the
composition $\lsplt \circ [E]: A\to \hf^1(K,E_2)$ can be described in
terms of Coleman integration.

Assume that we are given a curve
$X$ over $K$, zero divisors $y$ and $z$ with disjoint supports $\cy$ and
$\cz$ respectively. We assume that $y$ and $z$ are sums of
$K$-rational points. This implies that $A=\tA$ in the notation of
Section~\ref{sec:mixed}. This restriction will be eliminated in
Section~\ref{sec:end}.

In the situation we are considering, we have a natural isomorphism,
given by the exponential map of Bloch and Kato, 
\begin{equation}\label{exp}
  \ddr(V)/F^0 \xrightarrow{\exp} \hf^1(K,V)\;,
\end{equation}
where $\ddr$ is the de Rham functor of Fontaine, and
similarly with $V$ replaced by $E_2$. Recall that $\ddr(V)$ has a
filtration, which we used above, and is also an extension of scalars
from $K_0$, the maximal unramified extension of $\Q_p$ in $K$, of
$\dcr(V)$, which comes equipped with a semi-linear Frobenius. Thus,
$\ddr(V)$ comes with a natural Frobenius and the
same is true for $\ddr(E_2)$.

By the de Rham conjecture proved by
Faltings~\cite{Fal89}, we know that the de Rham functors applied to the local
representations are given by
$\ddr(V)=\hdr^1(X/K)$ and $\ddr(E_2)=\hdr^1(X; \cz/K)$, where in both
cases the cohomology is twisted by $1$ so that $F^0 \ddr = F^1 \hdr$.

An element $\beta\in \hf^1(K,E_2)$ defines a $K$-valued functional on $F^1
\hdr^1(X - \cz/K)$ as follows: By the exponential map \eqref{exp} we
may view $\beta$ as an element of $\hdr^1(X; \cz/K)/F^1$, which is
\poincare\ dual to $F^1 \hdr^1(X - \cz/K)$.

If $D\in A$, then $\lsplt\circ [E](D)\in \hf^1(K,E_2)$ gives such a
functional which we denote $[E]_D$. To describe this explicitly, note
that $F^1 \hdr^1(X - \cz/K)$ is the space of one-forms $\omega$ on $X$
with logarithmic singularities along $\cz$. Then we have
\begin{equation}
  \label{eq:blabla}
  [E]_D(\omega)=\omega\cup \exp^{-1}(\lsplt\circ [E](D)) 
\end{equation}
The following result is an extension of~\cite[Proposition 9.2]{Bes97}.
\begin{proposition}\label{relcol}
  For $D\in A$ we have $[E]_D(\omega)= \int_D\omega$.
\end{proposition}
To prove this result, we first recall the generalized Jacobian
$J$ of $X$ with respect to $\cz$. This can be described in two
different ways as
either the Picard group of line bundles of degree $0$ on $X$ with
trivialization
along $\cz$, or as divisors of $X$ with support disjoint from $\cz$
modulo divisors of functions that take the value $1$ on $\cz$. The maps
between the two descriptions takes a divisor $D$ to $\O(D)$ together
with the trivialization furnished by the canonical section $1$ and
going the other way take a line bundle trivialized along $\cz$ to the
divisor of a section which takes the value $1$ at all points of $\cz$
under the trivialization. Since we are assuming that $\cz$ splits over
$K$ there exists a short exact sequence
\begin{equation}\label{shortJ}
  0\to (\oplus_{x\in \cz} \gm)/\gm \to J\to J_X\to 0\;,
\end{equation}
where $J_X$ is the Jacobian of $X$.
This is given in the first description as follows: The left hand map
sends a collection
$(\alpha_x)\in \gm$ to the trivial line bundle with the
trivialization being the obvious one skewed by $\alpha_x$ at $x$. The
right hand map forgets the trivializations. In the second description
the right hand maps takes the class of the divisor $D$ to the class of
the same divisor (but under a stronger relation) and the left hand map
sends the collection $(\alpha_x)$ to the divisor of a function that
takes the value $\alpha_x$ at $x$.

The following is known:
\begin{lemma}\label{Freitag}
  the Galois representation $E_2=\het^1(\Xbar;\bsZ,1) $ is isomorphic to
   $V_p(J)$, the 
   $p$-adic Tate-module of $J$ tensored with $\Q_p$. Furthermore, the
   exact sequence \eqref{shortB} is induced by the short exact
   sequence \eqref{shortJ}.
\end{lemma}
\begin{proposition}\label{Raskind}
  Let $D$ be a zero-divisor on $X$. Then we have
  $[E](D)=\alpha_J([D]) \in H^1(K,V_p(J))$ where
  $\alpha_J$ is the Kummer map of $J$ and $[D]$ is the class of $D$ in
  $J(K)$.
\end{proposition}
\begin{proof}
This follows by the same argument as the one used in the proof
of~\cite[Appendix, Lemma]{Ras95}
\end{proof}
The following is an easy extension of the well known identification of
the tangent space to the Picard scheme (probably, also well known).
\begin{lemma}\label{tangent}
  There is a canonical isomorphism $T_0(J) \to \hdr^1(X;\cz /K)/ F^1$.
\end{lemma}
By~\cite[Example 3.10.1]{Blo-Kat90} there exist a commutative diagram
\newcommand{\Lie}{\operatorname{Lie}}
\begin{equation}\label{exkummer}
\xymatrix{
\hdr^1(X;\cz/K)/F^1=T_0(J) \ar[r]^(.65){\exp_J} \ar[rd]^{\exp}
&J(K)\otimes \Q_p
\ar[d]^{\alpha_J}\\
& \hf^1(K,E_2)
}
\end{equation}
Where $\exp_J$ is the exponential map of the $p$-adic Lie group $J(K)$
while $\exp$ is the Bloch-Kato exponential.
\begin{proposition}\label{fromSerre}
  Assume that $X$ has a $K$-rational point $P_0$ and let $g: X-\cz \to
  J$ be the map $P\mapsto P-P_0$. Then the map $g^\ast: \Omega^1(J) \to
  F^1 \hdr^1(X-\cz /K)$ is an isomorphism (recall that the space on
  the right is the
  space of forms with logarithmic singularities along
  $\cz$). Furthermore, this isomorphism fits in a commutative diagram
  \begin{equation*}
    \xymatrix{
      \Omega^1(J)\ar[r]^(.35){g^\ast}\ar[d]_{\operatorname{ev}_0} &
      F^1 \hdr^1(X-\cz /K)
      \ar[d]^{\text{\poincare\ duality}} \\
      T_0(J)^\ast \ar[r] & \hdr^1(X;\cz /K)^\ast
      }
  \end{equation*}
  where $\operatorname{ev}_0$ is evaluation at $0$ and the bottom map is the
  dual of the isomorphism of Lemma~\ref{tangent}.
\end{proposition}
\begin{proof}
  The first statement is~\cite[Proposition~V.5]{Ser59}. The second
  statement is again well known for the usual jacobian (see, e.g.,
  the proof of Proposition~2.2 in \cite{Mil86a}) and the
  proof is the same.
\end{proof}
\begin{proof}[Proof of Proposition~\ref{relcol}]
The short exact sequence \eqref{shortJ} induces a short exact sequence
of tangent spaces compatible with exponential maps. In addition, the
exponential map for $J_X$ is onto. Even more precisely, an integer
multiple of every class in
$J_X(K)$ is in the image of the exponential map, without tensoring
with $\Q_p$. This implies that any divisor $D$ on $X$, can be, after
multiplying by a sufficiently large integer, written as a sum of a
divisor, whose class, $[D]\in J(K)$, is in the image of $\exp_J$, and a
principal divisor, both disjoint from $\cz$. It is thus sufficient to
check the result for these two types of divisors.

We now identify $T_0(J)$ with $\hdr^1(X,\cz/K)$ and $\Omega^1(J)$ with
$F^1 \hdr^1(X-\cz /K)$. The duality between these spaces is given by
the cup product.

For the first type the proof is essentially the same as
the proof of~\cite[Proposition~9.2]{Bes97}: By assumption we have
$[D]= \exp_J(\eta)$ for some $\eta\in \hdr^1(X;\cz/K)/F^1$.
We have
\begin{align*}
  [E](D) &= \alpha_J([D]) \quad \text{by Proposition~\ref{Raskind}}\\
         &= \alpha_J( \exp_J(\eta))\\
         &= \exp(\eta) \quad \text{by \eqref{exkummer}.}
\end{align*}
In particular we have $[E](D)\in \hf^1(K,E_2)$ so that $\lsplt\circ
[E](D)= [E](D)=\exp(\eta)$ as well. It follows that
$[E]_D(\omega)=\omega\cup \eta$.
By Proposition~\ref{fromSerre} there
exists a holomorphic invariant differential
$\omega'$ on $J$ such that $g^\ast \omega'= \omega$. By
functoriality of the
$p$-adic integral (Proposition~\ref{funct}) we have $\int_D \omega=
\int_0^{[D]} \omega'$, where $[D]$ is the image of $D$ in $J(K)$.
Let $\log_J:J(K)\to T_0(J)=\hdr^1(X;\cz /K)$ be the
logarithm for the group $J(K)$ (depending on the choice of $\log_\uni$). 
By Corollary~\ref{funct1}, $\int_0^{[D]}
\omega'=\pair{\omega',\log_J([D])}=\omega\cup \log_J([D])$, where
$\pair{~,~}$ is the duality between the tangent space and holomorphic
forms on the jacobian.
But $\log_J([D])=\eta$ so the proof in this case is complete.

Consider now the case where $D$ is the divisor of a rational function
$f$. If $f$ take values
in $\O_K$ on $\cz$, then $D$ belongs to the image of $\exp_J$. Thus we
may assume that $\log(f)=0$ on $\cz$. From the description of the
maps in \eqref{shortJ} and Lemma~\ref{Freitag} it follows that
$[E](D)=j(\oplus_{x\in \cz} \alpha_K(f(x)))$. The definition of
$\lunip$ and the fact that $\log(f)=0$ on $\cz$ implies that
$1\otimes \lunip(\oplus_{x\in \cz} \alpha_K(f(x)))=0$. It follows from
Lemma~\ref{4.2} that $\lsplt\circ [E](D)=0$. Now we apply the theory of
the double index in the algebraic version of~\cite{Bes00} (in this
case it can also be deduced with a bit of effort from the rigid theory
of~\cite{Bes98b}). We take the sum of the double indices $\pair{\int
\omega,\log(f)}_x$ of the two $p$-adic integrals $\log(f)$ and
$\int(\omega)$ over all points $x$. A general result (in the rigid
case see~\cite[Corollary~4.11]{Bes98b}) says that the sum of these
indices for two functions,
one of which is a log, is always $0$. By definition, the sum of the
local indices is $\int_{(f)} \omega - \log f(\res(\omega))$ so this
last expression is $0$.
But by our assumptions $\log f(\res(\omega))=0$ so also $\int_{(f)}
\omega=0$ and the proof is complete.
\end{proof}

\section{End of the proof}
\label{sec:end}

To finish the description of the height pairing and complete the proof
we must describe the map
$w: \hf^1(K,E_2) \to \hf^1(K,B(1))$, giving a splitting of the first
line of \eqref{bigdiag}.
To describe this splitting is equivalent to describing a splitting 
$\hf^1(K,V)\to \hf^1(K,E_2)$. By the fact that the exponential
map~\eqref{exp} is an isomorphism for $V$ and $E_2$ this is the same
as describing a
section $\ddr(V)/F^0 \to \ddr(E_2)/F^0$ to the natural map in the
other direction.

Under the assumptions in~\cite{Nek93}, which are satisfied in our
case, there is a unique
Frobenius equivariant splitting of $\ddr(E_2)\to \ddr(V)$. We now
throw in the space $W=W_v \in \ddr(V)$, complementary to $F^0$. This
gives the required splitting,
\begin{equation}
  \label{eq:split}
  \ddr(V)/F^0 \xrightarrow{\text{along } W} \ddr(V)
  \xrightarrow{\text{Frob equivariant}} \ddr(E_2) \to \ddr(E_2)/F^0\;.
\end{equation}
This completes the description of the local height pairing.

In our situation we must describe the map
\begin{equation}
  \label{eq:split1}
  \hdr^1(X/K)/F^1 \xrightarrow{\text{along } W} \hdr^1(X/K)
  \xrightarrow{\text{Frob equivariant}} \hdr^1(X;\cz/K) \to
  \hdr^1(X;\cz/K)/F^1\;. 
\end{equation}
Dualizing, we have the splitting
\begin{equation}
  \label{eq:split2}
  F^1 \hdr^1(X-\cz/K)
  \to
 \hdr^1(X-\cz/K)
  \xrightarrow{\text{Frob equivariant}}  \hdr^1(X/K)
\rightarrow
  F^1 \hdr^1(X/K)
\end{equation}
Let $W'$ be the annihilator of $W$ with respect to the cup product. It
is easy to check that the rightmost map is the projection with respect
to the direct sum
decomposition $ \hdr^1(X/K)=W' \oplus  F^1 \hdr^1(X/K)$.

\begin{proposition}\label{getpsi}
  The composition of the two leftmost maps in \eqref{eq:split2} is
  induced by the map $\lag$ of \eqref{eq:lag}
\end{proposition}
\begin{proof}
From the description of $\lag$ in Proposition~\ref{proplag} this
is clear from the following lemma.
\end{proof}
\begin{lemma}
  Let $X_s$ be the closed fiber of $\XX$. Let $S=\{s_1,\ldots,s_n\}$ be
  a finite set of closed points of $X_s$ such that every point of $\cz$
  reduces to a point of $S$. Then, the map $\hdr^1(X- \cz/K)\to
  \hr^1(X_s-S/K)$ is Frobenius equivariant.
\end{lemma}
\begin{proof}
Note that the Frobenius structure on $\hdr^1(X- \cz/K)$ is
defined via its identification with $\ddr$ of \etale\ cohomology and is
not directly related to any cohomology on which Frobenius acts
naturally. However, as both de Rham and rigid cohomology have
\poincare\ duality (but not cristalline cohomology, which is badly
behaved for open schemes), we can dualize and must show that the
natural map $\hr^1(X_s;S/K)\to \hdr^1(X; \cz/K)$ is Frobenius
equivariant. Lift each $z_i\in \cz$ to a section $z_i:\spec(\O_K) \to \XX$,
and let $z_i^s: \spec(\kappa) \to X_s$ be the reduction. Then, the comparison
isomorphisms in a relative situation (see for example~\cite{Kis00}) imply
that the the Frobenius structure on $ \hdr^1(X; \cz/K)$ is induced via
the isomorphism $ \hdr^1(X; \cz/K)\isom \hcr^1(X_s; \cup_i
\spec(\kappa)/K_0)\otimes K$, where $\hcr^1(X_s; \cup_i
\spec(\kappa)/K_0)$ is the
cohomology relative to the maps $z_i^s$. As both $X_s$ and
$\spec(\kappa)$ are proper the relative crystalline cohomology is the
same as the relative rigid cohomology and the lemma follows easily.
\end{proof}

\begin{corollary}\label{lastcor}
  The induced splitting $B^\ast \otimes K \to F^1 \hdr^1(X-\cz/K)$ is
  given by $z\to \omega_{W'}(z)$, where $\omega_{W'}(z)$ is the form
  defined by Definition~\ref{omd}.
\end{corollary}
\begin{proof}
 The induced splitting is given by the following recipe: for $z\in B^\ast$
 take $\omega\in  F^1 \hdr^1(X-\cz/K)$ whose residue divisor is $z$,
 apply the map  \eqref{eq:split2} and subtract the result from
 $\omega$. The corollary follows immediately.
\end{proof} 

\begin{proof}[Proof of Theorem~\ref{locthm}]
We are still assuming that $\cy$ is split over
$K$. Let $y\in A$. By \eqref{bigfactor} We need to apply $z\in B^\ast$ to
$h_{E}(y) = (1\otimes \ell_v) \circ
u\circ w \circ \lsplt \circ [E](y)$. By
Proposition~\ref{relcol}, $\lsplt \circ [E](y)$ is a functional on
$\hdr^1(X -  \cz/K)$ defined by integrating a form on $y$. According
to Corollary~\ref{lastcor} this is mapped by $w$ to the functional on
$B^\ast \otimes  \hf^1(K,\Q_p(1))\isom B^\ast \otimes K$ obtained by
applying this functional to the form $\omega_{W'}(z)$. Finally, we apply to
the result $(1\otimes \ell_v) \circ u$, which by \eqref{whertvcomes}
is just $1\otimes t_v$. This gives
$z(h_E(y))=t_v(\int_{y}\omega_{W_v^\prime}(z))$. By the
proof of Proposition~5.2 in~\cite{Col-Gro89} we
have $\int_{y}\omega_{W_v^\prime}(z)=\int_{z}\omega_{W_v}(y)$, so
$z(h_E(y))=h_{\CG}^v(y,z)$ as required.

To finish the proof we must remove the assumption that $\cy$ is split
over $K$. For this it is sufficient to observe that the two local height
pairings behave in the same with respect to finite extensions. Indeed,
suppose that $K_1$ is a finite extension of $K$. Let $\ell_1$ be an
extension of $\ell_v$ to $K_1^\times$. Such an extension can be given
by $t_1\circ \log_v$ where $\log_v$ extends uniquely to $K_1^\times$
and $t_1:K_1\to \Q_p$ extends $t_v$. Let $W_1=W_v \cdot K_1$. Let $E$
be a mixed extension over $K$. Write
$h_{E/K_1}$ for the height
pairing corresponding to the mixed extension $E$ restricted to $K_1$,
taken with respect to $\ell_1$ and $W_1$. Then it is straightforward
to see that $h_{E/K_1}=h_E$. Thus, the local height pairing is
unchanged when we extend the field to $K_1$, provided that we choose
the data $(\ell_1,W_1)$ as above. If we now compare the Coleman-Gross
height pairing we see that the form $\omega_z$ is unchanged. Since the
log is also unchanged the integral $\int_y \omega_z$ is
unchanged. Finally, since this integral is in fact in $K$ and $t_1$
coincides with $t_v$ on $K$ we see that the entire local height
pairing is unchanged by our extension of scalars. This completes the
proof.
\end{proof}

\end{document}